\documentclass[12pt,a4paper]{article}
\usepackage[top=3.5cm, bottom=3.5cm, left=3cm, right=3cm]{geometry}
\usepackage[latin1]{inputenc}
\usepackage{csquotes}
\usepackage{amsmath}
\usepackage{amsthm}
\usepackage{amsfonts}
\usepackage{amssymb}
\usepackage{graphics}
\usepackage{float}
\usepackage[hang,flushmargin]{footmisc} 

\usepackage{amsmath,amsthm,amssymb,graphicx, multicol, array}
\usepackage{enumerate}
\usepackage{enumitem}
\usepackage{xcolor}

\usepackage[pagebackref,bookmarks,colorlinks,breaklinks]{hyperref}
\hypersetup{linkcolor=blue,citecolor=red,filecolor=blue,urlcolor=blue} 

\usepackage{tabto}

\numberwithin{equation}{section}

\usepackage{tikz}
\usepackage{pgfplots}

\usetikzlibrary{matrix}

\usepackage{mathrsfs}

\DeclareMathOperator{\lcm}{lcm}

\newtheorem{thm}{Theorem}[section]
\newtheorem{lem}{Lemma}[section]

\newtheorem{conj}{Conjecture}[section]

\newcommand{\Z}{\mathbb{Z}}
\newcommand{\Q}{\mathbb{Q}}

\newcommand{\C}{\mathbb{C}}
\newcommand{\F}{\mathbb{F}}

\usepackage{fancyhdr}
\title{Prime Values of the Euler Polynomial}
\date{}
\author{N. A. Carella}

\begin{document}

\thispagestyle{empty}
\date{}
\maketitle

\begin{abstract}
This note proposes an effective lower bound for the number of  primes in the quadratic progression $p=n^2+1 \leq x$ as $x \to \infty$.\let\thefootnote\relax\footnote{ \today \date{} 
	}
\end{abstract}


\section{Introduction} \label{S9090}\hypertarget{S9090}
As early as 1760, Euler was developing the theory of prime values of polynomials. In fact, Euler computed a very large table of the primes $p=n^2+1$, see {\color{red}\cite[p.\ 123]{EL1900}}. Likely, the prime values of polynomials was studied by other researchers before Euler. Later, circa 1910, Landau posed an updated question of the same problem about the primes values of the polynomial $n^2+1$. A fully developed conjecture, based on circle methods, was demonstrated about two decades later.

\begin{conj} \label{conj9090.112}\hypertarget{conj9090.112}{\normalfont ({\color{red}\cite[Conjecture E, p.\,46]{HL1923}})} Let $x \gg 1$ be a large number, let $\Lambda(n)$ be the vonMangoldt function, and let $\chi(n)=(n\mid p)$ be the quadratic symbol modulo $p$. Then
\begin{equation} \label{eq9090.100}
\sum_{n\leq x} \Lambda \left( n^2+1 \right )
=a_2x+O\left (\frac{ x}{\log x }\right ),
\end{equation}
where the density constant
\begin{equation} \label{eq9090.120}
a_2=\prod_{p \geq 3}\left( 1-\frac{\chi(-1)}{p-1} \right )=1.37281346 \ldots .
\end{equation}
\end{conj}

The claim in Conjecture E is specific to the polynomial $f(t)=t^2+1$ and the next conjecture F in \cite{HL1923} covers the general admissible quadratic polynomial $f(t)=at^2+bt+c\in \Z[t]$. An \textit{admissible} polynomial $f(t)\in \Z[t]$ has a relatively prime image $\gcd(f(\Z))=1$, see {\color{red}\cite[p.\;385]{FI2010}}, \cite{FM1988} and similar sources for detailed explanations. A survey of the subsequent developments appears in {\color{red}\cite[p.\ 342]{NW2000}}, {\color{red}\cite[Section 19]{PJ2009}}, \cite{BZ2007}, and similar references. Other related problems are discussed in \cite{BH1962}, \cite{BZ2007}, \cite{GL2010}, \cite{FI2010}, {\color{red}\cite[p. 343]{NW2000}}, {\color{red}\cite[p. 405]{RP1996}}, \cite{PJ2009}. \\

There are many approaches and techniques toward the proof of \hyperlink{conj9090.112}{Conjecture} \ref{conj9090.112} on the prime values of the admissible irreducible polynomial $n^2+1$ and other admissible and irreducible polynomials. Two of the leading techniques are essentially approximations to the actual problem: 
\begin{itemize}
\item Determining the \textit{largest prime factor }of the product $\prod_{n\leq x}(n^2+1)$. \\
\item Determining the \textit{smallest coefficient} $m$ in the representation $p=n^2+m^2$. \\
\end{itemize}

There is a large literature dedicated to the determination of the largest prime factors of admissible irreducible polynomials. The analysis of these results are based on the Chebyshev-Hooley method, which is a difficult analytical technique, see \cite{HC1967}, \cite{DI1982}, \cite{MJ2023}, et alii.  In addition, the properties of the primitive divisors of pure polynomials such as  $f(t)=t^2+c$ and some related topics are studied in \cite{EH2008}, \cite{HG2008}, \cite{HG2024}, et alii. An abridged table of the known results focused on the polynomial $f(t)=t^2+1$ are listed in \autoref{table9090-A}. Here, the symbol $P(n)\geq1$ denotes the largest prime factor of an integer $n\geq1$.

\vskip .1 in 
\begin{table}[H]
	\setlength{\tabcolsep}{0.5cm}
	\renewcommand{\arraystretch}{1.750092}
	\setlength{\arrayrulewidth}{0.80pt}
	\centering
	\begin{tabular}{l|l}
		
		Estimate of the largest prime divisor	& Author and reference \\
		\hline
		$P(n^2+1)/n\to \infty $ as $ n\to\infty$	&  Chebyshev, see \cite{HC1967}.\\
		\hline
		$P(n^2+1)\geq n(\log n)^{\varepsilon}$ as $ n\to\infty$	& Nagell, see \cite{NT1921}. \\
		\hline
		$P(n^2+1)\geq n(\log n)^{c\log\log\log n}$ as $ n\to\infty$		&  Erdos, see \cite{EP1952}.\\
		\hline
		$P(n^2+1)\geq n^{11/10}$ as $ n\to\infty$		&  Hooley, see \cite{HC1967}.\\
		\hline
		$P(n^2+1)\geq n^{1.202}$ as $ n\to\infty$		&  Iwaniec, see \cite{DI1982}.\\
		\hline
		$P(n^2+1)\geq n^{1.279}$ as $ n\to\infty$		&  Merikoski, see \cite{MJ2023}.\\
	\end{tabular}  
	\vskip .1 in		
	\caption{Estimated magnitude of prime divisor of $n^2+1$}
	\label{table9090-A}
\end{table}

Likewise, there is a large literature dedicated to the determination of smallest coefficient in the representation $p=n^2+m^2$. In terms of complex gaussian numbers $\pi=n+im\in\Z[i]$, the rational integer $n^2+m^2=\pi \overline {\pi}$ is a prime number infinitely often if and only if the complex line $\Im  m \;\pi=1\subset \C$ contains infinitely many gaussian primes $\pi=n+i$. Moreover, conditional on the generalized RH, there are a few proofs for the existence of infinitely many gaussian primes $\pi=a+ib\in \Z[i]$ on the complex strip $0\leq \Im  m \;\pi\ll \log p\subset \C$, confer {\color{red}\cite[Theorem III]{AN1952}}, et alii.  Unconditionally, there are various results regarding the existence of infinitely many gaussian primes $\pi=a+ib\in \Z[i]$ on the complex strip $0\leq \Im  m \;\pi\ll p^{\alpha}\subset \C$, where $\alpha<1/2$. More generally, the proof for quadratic forms $p=Q((m,n)=am^2+bmn+cm^2$, with $m\ll p^{\alpha}$, and a survey of previous results appears in \cite{CM1993}. These results are equivalent to the results for complex sectors $\alpha\leq \arg \pi \leq \beta\subset \C$. Some topics in the advanced algebraic and analytic theory of the gaussian number field $\Q(i)$ and primes as sums of two squares, such as the primes $p=n^2+1$, are treated in many papers, see \cite{FI2022}, {\color{red}\cite[Chapter 12]{HG2007}} et alii. An abridged table of the progress is listed on \autoref{table9090-B}.  
\vskip .1 in 
\begin{table}[H]
	\setlength{\tabcolsep}{0.5cm}
	\renewcommand{\arraystretch}{1.750092}
	\setlength{\arrayrulewidth}{0.80pt}
	\centering
	\begin{tabular}{l|l|l}
		Smallest coefficient $m$	& Criterion&Authors and references \\
		\hline
		$m\ll\log  p$	&Conditional  &Ankeny 1952, \cite{AN1952}\\		
		\hline
		$m\ll p^{1/4}$	& Unconditional& Simple exercise.\\
		\hline
		$m\ll p^{0.1631}$	&Unconditional  &Coleman 1993, \cite{CM1993}\\
		
		\hline
		$m\ll p^{0.119}$	&Unconditional  &Harman-Lewis 2001, \cite{HL2001}\\
	\end{tabular}
	\vskip .1 in		
	\caption{The prime $p=n^2+m^2$ and the smallest $m$}
	\label{table9090-B}
\end{table}

Some partial results are proved in \cite{GM2000}, \cite{MA2009}, \cite{BZ2007}, \cite{DI1982}, \cite{IH1978}, \cite{LR2012}, et alii. The smallest expected prime value of certain admissible polynomials $x^d+a\in\Z[x]$ is studied in \cite{MK1986}. Assuming the Elliott-Halberstam conjecture, there is a discussion in {\color{red}\cite[p.\ 5]{BZ2007}} concerning the existence of infinitely many primes of the form $p=an^2+1$, with $a =O(p^{\varepsilon})$, and $\varepsilon>0$. Extensive discussion on this later conjecture and the level of distribution of the moduli are given in {\color{red}\cite[p.\ 406]{FI2010}}. The results for the associated least common multiple problem $\log \lcm [f(1)f(2)\cdots f(n)]$ appears in \cite{CJ2011}, and the recent literature. This note proposes the following partial result.

\begin{thm} \label{thm9090.119}\hypertarget{thm9090.119} Let $x \geq 1$ be a large number. Then, 
\begin{equation} \label{eq909.19}
 \sum_{\sqrt{x}< n\leq \sqrt{2x}} \Lambda(n^2+1)\gg \sqrt{x}\left( 1+O\left(e^{-c\sqrt{\log x}} \right) \right).
\end{equation}
\end{thm}

This result seems to resolve the asymptotic part of the polynomial prime values problem for the admissible polynomial $f(x)=x^2+1$. The proof, based on standard elementary methods in number theory, appears in \hyperlink{S2099}{Section} \ref{S2099}. The supporting materials are developed in \hyperlink{S2070}{Section} \ref{S2070}, \hyperlink{S2099MT}{Section} \ref{S2099MT} and \hyperlink{S2099ET}{Section} \ref{S2099ET}. \\

Some heuristics, and discussions on the difficulty and complexity of estimating or computing the constant appear in \cite{BH1962}, {\color{red}\cite[Section 3.3]{FG2018}}, \cite{RI2015}, and other papers on the Bateman-Horn Conjecture.

\section{Linear to Quadratic Identity} \label{S2070}\hypertarget{S2070}
Among the representations of the characteristic function of square integers $n \in [1,x]$ there is the double finite sum
\begin{equation}\label{eq2070.130d}
\mathcal{I}_2(n)=\frac{1}{N}\sum_{m \leq \sqrt{x},}\sum_{0\leq u<N} e^{i2\pi \left(m^2-n \right)u/N } =
\begin{cases}
1&\text{  if }  n=m^2,\\
0&\text{  if }  n\ne m^2,\\
\end{cases} 
\end{equation}
where $N>1$ is a parameter and $x>1$ is a real number. The weighted characteristic function of prime number, better known as the vonMangoldt function, is defined by

\begin{equation}\label{eq2070.130f}
\Lambda(n)=
	\begin{cases}
		\log p &\text{  if }  n=p^v,\\
		0&\text{  if }  n\ne p^v,\\
	\end{cases} 
\end{equation}
where $n=p^v$ is a prime power.\\

Merging these two characteristic functions produces the
\textit{quadratic to linear identity}.

\begin{lem}\label{lem2070.130}\hypertarget{lem2070.130} If $n\geq$ is an integer, then
\begin{equation}\label{eq2070.130i}
\mathcal{I}_2(n)\Lambda(n+1) =
\begin{cases}
\Lambda(m^2+1)&\text{  if }  n=m^2 \text{ is a square,}\\
0&\text{  if }  n\ne m^2 \text{ is not a square.}\\
\end{cases} 
\end{equation}
\end{lem}  

The application of the characteristic function \eqref{eq2070.130d} presented here uses the choice of parameter $N\approx x$ or $N\sim x$ and $N$ prime. For a fixed $n\in[1,x]$, this ensures that the congruence equation 
\begin{equation} \label{eq2070.100b}
	z^2-n\equiv 0\bmod N ,
\end{equation}
where $z\leq \sqrt{x}$, and the integer equation 
\begin{equation} \label{eq2070.100c}
	z^2-n= 0
\end{equation} 
have exactly the same integer solution as specified in \eqref{eq2070.130d}. The matched  parameters ensures that the double exponential sum \eqref{eq2070.130d} works as a characteristic function of square integers $n\leq x$ as specified. Mismatched parameters produce \textit{quasi characteristic functions}. For example,  
\begin{equation}\label{eq2070.130g}
\frac{1}{N}\sum_{m \leq x^{3/4},}\sum_{0\leq u<N} e^{i2\pi \left(m^2-n \right)u/N }
\end{equation}
is not a characteristic function of the squared integers for large $N\approx x$.

\section{Fibers and Multiplicities} \label{S2095SS} \hypertarget{S2099SS}
The cardinalities and multiplicities of the fibers occurring in the estimate of the error term are computed in this section. The notation $[x]=x-{x}$ denotes the integer value  of the quantities $x>0$.
\begin{lem}  \label{lem2099SS.300S}\hypertarget{lem2099SS.300S} Let $N$ be a prime and define the maps
	\begin{equation}\label{eq2099SS.300-m}
		\alpha(m,n)\equiv (m^{2}-n)\bmod N\quad \text{ and } \quad 
		\beta(u,v)\equiv uv\bmod N.
	\end{equation}	
	If $\Lambda(n^2+1)=0$ for all $n\in (N,2N)$, then the fibers $\alpha^{-1}(r)$ and $\beta^{-1}(r)$ of an element $0\ne r\in \F_N$
	have the cardinalities 
	\begin{equation}\label{eq2099SS.300-f}
		\#	\alpha^{-1}(r)\leq [\sqrt{2N}]-[\sqrt{N}]-1\quad \text{ and }\quad \#\beta^{-1}(r)=	[\sqrt{N}]
	\end{equation}
	respectively.
\end{lem}
\begin{proof}[\textbf{Proof}] Given a fixed $m\in \mathscr{M}=\{\sqrt{N}<m\leq \sqrt{2N}\}$, the map 
	\begin{equation}\label{eq9900Q.300-m1}						  		
		\alpha: \mathscr{M}\times[1,N-1]\longrightarrow\F_N\quad  \text{ defined by }\quad  \alpha(m,n)\equiv (m^{2}-n)\bmod N,
	\end{equation}
	where $n\ne m^2$,	is one-to-one but not surjective since the value 
	\begin{equation}\label{eq9900Q.300-m2}						  		
		\alpha(m,n)= (m^{2}-n)\equiv 0 \bmod N
	\end{equation} is missing. This follows from the fact that the map $n\longrightarrow m^2-n \bmod N$ is a permutation the nonezero elements of the finite field $\F_N$,
	see {\color{red}\cite[Chapter 7]{LN1997}} for more details on permutation functions of finite fields. Thus, as $(m,n)\in \mathscr{M}\times[1,N-1]$ varies, a value $0\ne r=\alpha(m,n)\in\F_N$ is repeated at most $[\sqrt{2N}]-[\sqrt{N}]-1$ times since the hypothesis implies that $m^2-n\ne0$ for all $(m,n)$. This verifies that the cardinality of the fiber is
	\begin{eqnarray}\label{eq9900Q.300-f1}
		\#	\alpha^{-1}(r)&=&	\#\{(m,n):r\equiv (m^{2}-n)\bmod N:1\leq m\leq \sqrt{N} \text{ and } n<N/2\}\nonumber\\[.3cm]
		&=& [\sqrt{2N}]-[\sqrt{N}]-1.
	\end{eqnarray}		
	Similarly, given a fixed $u\in [1,\sqrt{N}]$, the map 
	\begin{equation}\label{eq9900Q.300-m3}
		\beta:[1,\sqrt{N}]\times [1,N-1]\longrightarrow\F_N\quad  \text{ defined by }\quad  \beta(u,v)\equiv uv\bmod N,
	\end{equation}
	is one-to-one. Here the map $v\longrightarrow uv \bmod N$ permutes the nonzero elements of the finite field $\F_N$. Thus, as $(u,v)\in [1,\sqrt{N}]\times [1,N-1]$ varies, each value $0\ne r=\beta(u,v)\in\F_N$ is repeated exactly $[\sqrt{N}]$ times. This verifies that the cardinality of the fiber is
	\begin{eqnarray}\label{eq9900Q.300-f2}
		\#	\beta^{-1}(r)&=&	\#\{(u,v):r\equiv uv\bmod N:1\leq u\leq \sqrt{x} \text{ and }1\leq v< N\}\nonumber\\[.3cm]
		&=&[\sqrt{x}].
	\end{eqnarray}
	
	Now each value $r=\alpha(m,n)\ne0$ (of multiplicity at most $[\sqrt{2N}]-[\sqrt{N}]-1$ in $	\alpha^{-1}(r)$), is matched to $r=\alpha(m,n)=\beta(u,v)$ for some $(u,v)$, (of multiplicity exactly $[\sqrt{N}]$ in $	\beta^{-1}(r)$). Comparing \eqref{eq9900Q.300-f1} and \eqref{eq9900Q.300-f2} proves that
	\begin{equation}\label{eq9900Q.300-m4}
		\# \alpha^{-1}(r)=[\sqrt{2N}]-[\sqrt{N}]-1<\# \beta^{-1}(r)=[\sqrt{N}].
	\end{equation}
\end{proof}

\section{Estimate for the Main Term} \label{S2099MT}\hypertarget{S2099MT}
In this analysis the relationship $N\sim x$ between the real number $x$ and the prime $N$ is explicated in \hyperlink{S2070}{Section} \ref{S2070}. 
\begin{lem} \label{lem2099MT.100A}\hypertarget{lem2099MT.100A} If $x \geq 1$ is a large number and $N\sim x$ is a large prime, then
\begin{equation} \label{eq2099.100d}
\frac{1}{N}\sum_{\substack{x< n\leq 2x\\n\ne m^2}} \Lambda(n+1) \sum_{\sqrt{x}< m\leq \sqrt{2x}}1\gg \sqrt{x}\left( 1+O\left(e^{-c\sqrt{\log x}} \right) \right),\nonumber
\end{equation}
where $c>0$ is a constant.	
\end{lem}

\begin{proof}[\textbf{Proof}] 
Set the parameter $N\sim x$ as a prime greater than $x$. By hypothesis the subsum
\begin{equation} \label{eq2099.100f}
\sum_{\sqrt{x}< m\leq \sqrt{2x}} \Lambda(m^2+1)=\sum_{x<m^2=n\leq 2x} \Lambda(n+1)=0
\end{equation} 
does not contribute to the main term since there are no primes of the form $p=n^2+1\in (x,2x]$. Thus, the main term can be rewritten and estimated as 
\begin{eqnarray} \label{eq2099.100h}
M(x)&=&\frac{1}{N}\sum_{x< n\leq 2x} \Lambda(n+1)\sum_{\sqrt{x}< m\leq \sqrt{2x}}1\\	[.3cm]
	&=&\frac{(\sqrt{2}-1)\sqrt{x}+O(1)}{N}\sum_{x< n\leq 2x}  \Lambda(n+1)  \nonumber\\[.3cm]
	&=&\frac{(\sqrt{2}-1)\sqrt{x}+O(1)}{N}\left( x+O\left(xe^{-c\sqrt{\log x}} \right) \right) \nonumber\\[.3cm]
		&\gg&\sqrt{x}\left( 1+O\left(e^{-c\sqrt{\log x}} \right) \right) \nonumber,
\end{eqnarray}
where $c>0$ is a constant. The third line in \eqref{eq2099.100h} follows from the prime number theorem, see {\color{red}\cite[27.12.E6]{DLMF}}, \cite{EL1985}, \cite{IK2004}, {\color{red}\cite[Theorem 6.9]{MV2007}}. Lastly, the last line follows from $N\sim x$.
\end{proof}

\section{Estimate for the Error Term} \label{S2099ET} \hypertarget{S2099ET}
The analysis of the estimate computed here is based on standard techniques in analytic number theory.
\begin{lem} \label{lem2099ET.100B}\hypertarget{lem2099ET.100B} Let $x \geq 1$ be a large number and let $N\sim x$ be a large prime. Suppose there are no primes $p=n^2+1$ in the short interval $(x,2x]$. Then
\begin{equation} \label{eq2099.400d}
E(x)=\sum_{\substack{x< n\leq 2x\\n\ne m^2}} \Lambda(n+1) \left (\frac{1}{N}\sum_{\sqrt{x}< m\leq \sqrt{2x},}\sum_{1\leq u<N} e^{i2\pi \left(m^2-n \right)u/N }\right) = O\left ( (\log x)^{3}\right )\nonumber.
\end{equation}
\end{lem}
\begin{proof}[\textbf{Proof}] Partitioning the inner finite sum in the following way.
\begin{eqnarray} \label{eq2099.400g}
	E(x)	&=&\sum_{\substack{x< n\leq 2x\\n\ne m^2}} \Lambda(n+1) \left (\frac{1}{N}\sum_{\sqrt{x}< m\leq \sqrt{2x},}\sum_{1\leq u<N} e^{i2\pi \left(m^2-n \right)u/N }\right) \\[.3cm]
	&=&\frac{1}{N}\sum_{\substack{x< n\leq 2x\\n\ne m^2}} \Lambda(n+1) \sum_{\sqrt{x}< m\leq \sqrt{2x},}\left (\sum_{1\leq u<N/2} e^{i2\pi \left(m^2-n \right)u/N }+\sum_{N/2\leq u<N} e^{i2\pi \left(m^2-n \right)u/N }\right)  \nonumber\\[.3cm]
	&=&E_0(x) \; + \; E_1(x)\nonumber.
\end{eqnarray}
Summing the inner sum in the first subsum yields
\begin{eqnarray} \label{eq2099.400i}
	E_0(x)	&=&\frac{1}{N}\sum_{\substack{x< n\leq 2x\\n\ne m^2}} \Lambda(n+1)\sum_{\sqrt{x}< m\leq \sqrt{2x},}\sum_{1\leq u<N/2} e^{i2\pi \left(m^2-n \right)u/N }\\[.3cm]
	&=& \frac{1}{N}\sum_{\substack{x< n\leq 2x\\n\ne m^2}} \Lambda(n+1) \sum_{\sqrt{x}< m\leq \sqrt{2x},} \frac{e^{i2\pi (\frac{m^{2}-n}{N})(\frac{N+1}{2})}-1}{e^{i2\pi \frac{(m^{2}-n)}{N}}-1} \nonumber\\[.3cm]
	&\leq &\frac{1}{N}\sum_{\substack{x< n\leq 2x\\n\ne m^2}} \Lambda(n+1) \sum_{\sqrt{x}< m\leq \sqrt{2x}}\left(\frac{2}{|\sin\pi(m^2-n)/N|} \right)	\nonumber.
\end{eqnarray}
The last inequality in \eqref{eq2099.400i} follows from the hypothesis $m^2-n\ne0$ and taking absolute value.\\

Replace $N\sim x$ and apply \hyperlink{lem2099SS.300S}{Lemma} \ref{lem2099SS.300S} to obtain the next inequality 
\begin{eqnarray} \label{eq2099.400t}
	E_{0}(x)&\leq &\frac{1}{N}\sum_{\substack{N< n\leq 2N\\n\ne m^2}} \Lambda(n+1) \sum_{\sqrt{N}< m\leq \sqrt{2N}}\left(\frac{2}{|\sin\pi(m^2-n)/N|} \right)	\\[.3cm]
	&\ll&   \frac{\log 2N}{N}\sum_{a\leq \sqrt{2N},}  \sum_{b \leq N}\left(\frac{1}{|\sin\pi ab/N|} \right)	\nonumber ,
\end{eqnarray} 
where $|\sin\pi ab/N|\ne0$ since $N\nmid ab$. Observe that as the index in the first double sum in \eqref{eq2099.400t} ranges over $(m,n)\in \mathscr{M}\times (N,2N]$, it ranges over $\sqrt{2N}-\sqrt{N}-1$ copies of the set $(N,2N]-\{m^2\}$, where $\mathscr{M}=\{\sqrt{N}+1,\sqrt{N}+2,\ldots, \sqrt{2N}\}$. In contrast, as the index in the last double sum in \eqref{eq2099.400t} ranges over $(a,b)\in [1,\sqrt{2N}]\times [1,N]$ it ranges over $\sqrt{2N}$ copies of the set $[1,N]$. Thus, \hyperlink{lem2099SS.300S}{Lemma} \ref{lem2099SS.300S} fully justifies the inequality \eqref{eq2099.400t}. In light of this information, continuing the calculation of the estimate as in {\color{red}\cite[p.\;136]{DH2000}}, yields 

\begin{eqnarray} \label{eq2099.400u}
	E_0(x)	
	&\ll &\frac{\log N}{N}\sum_{a\leq \sqrt{2N},}  \sum_{b \leq N}\left(\frac{1}{|\sin\pi ab/N|} \right)	\\[.3cm]	 
	&\ll &\frac{\log N}{N}\sum_{a\leq \sqrt{2N},}  \sum_{b \leq N}\left(\frac{N}{\pi ab} \right)	\nonumber\\	 [.3cm] 
	&\ll &(\log N)\sum_{a\leq \sqrt{2N}}\frac{1}{a}  \sum_{b \leq N}\frac{1}{b}	\nonumber\\	[.3cm]  
	&\ll &(\log N)^3	\nonumber\\[.3cm]
		&\ll &(\log x)^3	\nonumber.
\end{eqnarray}
Similarly, the second subsum has the upper bound
\begin{eqnarray} \label{eq2099.400v}
	E_{1}(x)&=& \frac{1}{N}\sum_{\substack{x< n\leq 2x\\n\ne m^2}} \Lambda(n+1) \sum_{m \leq \sqrt{x},}\sum_{N/2\leq u<N} e^{i2\pi \left(m^2-n \right)u/N } \\[.3cm]
	&=&    \frac{1}{N}\sum_{\substack{x< n\leq 2x\\n\ne m^2}} \Lambda(n+1) \sum_{m \leq \sqrt{x},} \frac{1-e^{i2\pi (\frac{m^{2}-n}{N})(\frac{N+1}{2})}}{e^{i2\pi \frac{(m^{2}-n)}{N}}-1} \nonumber\\[.3cm]
	&\leq&   \frac{1}{N}\sum_{\substack{x< n\leq 2x\\n\ne m^2}} \Lambda(n+1) \sum_{m \leq \sqrt{x}}\left(\frac{2}{|\sin\pi(m^2-n)/N|} \right) \nonumber\\[.3cm]
	&\ll&  (\log x)^3\nonumber.
\end{eqnarray}
Adding \eqref{eq2099.400u} and \eqref{eq2099.400v} and substituting \( x=(\log p)^{1+\varepsilon}\) yield
\begin{eqnarray} \label{eq2099.300x}
	E(x)&=& E_{0}(x)\;+\;E_{1}(x)  \\[.3cm]
	&\ll& (\log x)^3\nonumber.
\end{eqnarray}
This completes the verification.
\end{proof}

\section{Main Result} \label{S2099}\hypertarget{S2099}
The proof of the main result is based on the quadratic to linear identity described \hyperlink{S2070}{Section} \ref{S2070}.
\begin{proof}[\textbf{Proof}] (\hyperlink{thm9090.119}{Theorem} \ref{thm9090.119}) Take a sufficiently large number $x$ and suppose there are no primes $p=n^2+1\geq x$. Equivalently  $\Lambda(n^2+1)=0$ for all $n\geq x^{1/2}$. Summing the quadratic to linear identity in \hyperlink{lem2070.130}{Lemma} \ref{lem2070.130} over the short interval $(x,2x]$ leads to null sum
	\begin{equation} \label{eq2099.134d}
		\sum_{\sqrt{x}< m\leq \sqrt{2x}} \Lambda(m^2+1)= \sum_{\substack{x< n\leq 2x\\n\ne m^2}}\mathcal{I}_2(n)\Lambda(n+1)=0.
	\end{equation}
	This is a quantitative form of the hypothesis``\textit{there are no primes $p=n^2+1$ on the short interval} $(x,2x]$".\\
	
	Substituting the complete formula for the indicator function in \eqref{eq2070.130d} and expanding it yield 
	\begin{eqnarray} \label{eq2099.134f}
		\sum_{\sqrt{x}< m\leq \sqrt{2x}} \Lambda(m^2+1)
		&=& \sum_{\substack{x< n\leq 2x\\n\ne m^2}} \Lambda(n+1) \left (\frac{1}{N} \sum_{\sqrt{x}< m\leq \sqrt{2x},}\sum_{0\leq u<N} e^{i2\pi \left(m^2-n \right)u/N }\right) \nonumber\\[.3cm]
		&=&\frac{1}{N}\sum_{\substack{x< n\leq 2x\\n\ne m^2}} \Lambda(n+1) \sum_{\sqrt{x}< m\leq \sqrt{2x}}1 \nonumber\\[.3cm]
		&&\hskip .0005 in +\sum_{\substack{x< n\leq 2x\\n\ne m^2}}\Lambda(n+1) \left (\frac{1}{N} \sum_{\sqrt{x}< m\leq \sqrt{2x},}\sum_{1\leq u<N} e^{i2\pi \left(m^2-n \right)u/N }\right)   \nonumber\\[.3cm]
		&=&M(x) \; + \; E(x).
	\end{eqnarray}
	The main term $M(x)$, which is associated with the value $u=0$, is computed in \hyperlink{lem2099MT.100A}{Lemma} \ref{lem2099MT.100A}. The error term, which is associated with the values $u=1,\ldots,N-1$ is estimated in \hyperlink{lem2099ET.100B}{Lemma} \ref{lem2099ET.100B}. Substituting these estimates yield
	\begin{eqnarray} \label{eq2099.134i}
		\sum_{\sqrt{x}< m\leq \sqrt{2x}} \Lambda(m^2+1)
		&=&M(x) \; + \; E(x)\\[.3cm]
		&\gg& \sqrt{x}\left( 1+O\left(e^{-c\sqrt{\log x}} \right) \right)+O\left( (\log x)^3\right) \nonumber\\[.3cm]
		&\gg& \sqrt{x}\left( 1+O\left(e^{-c\sqrt{\log x}} \right) \right) \nonumber\\[.3cm]
		&\gg& 1 \nonumber
	\end{eqnarray}
	for all large real numbers $x$, where $c>0$ is a constant. But this contradicts the hypothesis \eqref{eq2099.134d}. Therefore, the short interval $[x,2x]$ contains a prime $p=n^2+1$ for all large real numbers $x$.  
	Quod erat inveniendum. \end{proof}

The earliest numerical data seems to be the Euler table in {\color{red}\cite[p.\ 123]{EL1900}}. Some other numerical data and experiments are reported in {\color{red}\cite[p.\ 50]{HL1923}}, \cite{SD1959}, and \cite{RI2015}. An instructive numerical experiment for the polynomial $f(x)=x^2+1$ is conducted in {\color{red}\cite[Section 3.3]{FG2018}}.



\end{document}